\documentclass[10pt]{article}
\usepackage{latexsym}
\usepackage{amsbsy}
\usepackage{amssymb}
\usepackage{amsmath}
\usepackage[latin1]{inputenc}
\usepackage{a4wide}
 \newcommand{\eps}{\varepsilon}
 \newcommand{\R}{\mathbb R}

 \newcommand{\gs}{\ensuremath{{\mathcal G}} }

 \newcommand{\esm}{\ensuremath{{\mathcal E}_M} }

 \newcommand{\comp}{\subset\subset}
 \newcommand{\cinfty}{{\cal C}^\infty}

 \newcommand{\vphi}{\varphi}

 \newcommand{\D}{{\cal D}}

 \newcommand{\G}{{\cal G}}

 \newtheorem{thr}{\hspace*{-3mm} \bf}[section]
 \newcommand{\ethr}{\end{thr}}

 \newcommand{\bd}{\begin{thr} {\bf Definition. }}
 \newcommand{\ed}{\end{thr}}

 \newcommand{\pr}{\noindent{\bf Proof. }}
 
 \newcommand{\ep}{\hspace*{\fill}$\Box$\ms}

 \newcommand{\bthm}{\begin{thr} {\bf Theorem. }}
 \newcommand{\bt}{\begin{thr} {\bf Theorem. }}
 \newcommand{\ethm}{\end{thr}}

 \newcommand{\bp}{\begin{thr} {\bf Proposition. }}
 \newcommand{\bc}{\begin{thr} {\bf Corollary. }}
 \newcommand{\blem}{\begin{thr} {\bf Lemma. }}
 \newcommand{\bex}{\begin{thr} {\bf Example. }\rm}
 \newcommand{\bexs}{\begin{thr} {\bf Examples. }\rm}
 \newcommand{\brm}{\begin{thr} {\bf Remark. }\rm}

 \newcommand{\beast}{\begin{eqnarray*}}
 \newcommand{\eeast}{\end{eqnarray*}}
 
 \newcommand{\brem}{\begin{thr} {\bf Remark. }\rm}
 \newcommand{\ethi}{\end{thr}}

 \newcommand{\ms}{\medskip\\}
 \newcommand{\nn}{\nonumber}
 \newcommand{\beq}{ \begin{equation} }
 \newcommand{\eeq}{\end{equation} }
 \newcommand{\bea}{\begin{eqnarray}}
 \newcommand{\eea}{\end{eqnarray}}
 \newcommand{\beas}{\begin{eqnarray*}}
 \newcommand{\eeas}{\end{eqnarray*}}
 \newcommand{\beqs}{\begin{equation*}}
 \newcommand{\eeqs}{\end{equation*}}

 \newcommand{\ben}{\begin{enumerate}}
 \newcommand{\een}{\end{enumerate}}
 \newcommand{\ba}{\begin{array}}
 \newcommand{\ea}{\end{array}}
 \newcommand{\cd}{{\cal D}}

 \newcommand{\cG}{{\cal G}}

\begin{document}



\title{Group invariants in algebras of generalized functions}
\author{Sanja Konjik \footnote{Faculty of Agriculture, University of Novi Sad, Trg Dositeja Obradovi\' ca 8, 21000 Novi Sad, SCG.
Electronic mail: kinjoki@neobee.net}\\
         Michael Kunzinger \footnote{Faculty of Mathematics, University of Vienna, Nordbergstr.\ 15, A-1090 Wien, Austria,
         Electronic mail: michael.kunzinger@univie.ac.at}\\
       }

\date{}
\maketitle

\begin{abstract}
We study invariance properties of Colombeau generalized functions
under  actions of smooth Lie transformation groups. Several
characterization results analogous to the smooth setting are derived
and applications to generalized rotational invariance are given.

\vskip5pt
\noindent
{\bf Mathematics Subject Classification (2000):} Primary: 46F30;
secondary: 58E40,46T30, 35A30
\vskip5pt
\noindent
{\bf Keywords:} Lie group actions, Colombeau generalized functions, invariants
\end{abstract}

\section{Introduction} \label{introduction}
Extensions of classical Lie group analysis (\cite{Olv,Olv2,BK}) to spaces
of generalized functions,
first in the distributional setting and later in nonlinear theories of
generalized functions have been developed by various authors, starting
as early as the 1950ies (\cite{Me, S, Zie, ga, bi, OR, RR, RW, symm, DKP,
ObCont, MoRot}). The main focus of the extension to Colombeau algebras
of generalized functions so far has been the study of symmetry groups
of differential equations involving singularities. In \cite{ObCont, MoRot}
(see also \cite{book}, ch.\ 4), however, M.\ Oberguggenberger initiated
the study of group invariant Colombeau generalized functions. He studied
invariance under translations and rotations and gave some applications
to the determination of group invariant fundamental solutions. Such
questions have been taken up in \cite{PSV, KoK} and will also be one
of the main themes of this contribution. Our aim is to analyze the
action of smooth Lie transformation groups on elements of Colombeau
algebras and to derive characterizations of invariants under such
transformations which parallel the analogous constructions in the
smooth setting.

The theory underlying our investigations is global analysis in Colombeau
algebras, as presented in \cite{book} (see \cite{horvath} for a
recent survey). In particular, our study of one-parameter transformation
groups is based on the theory of flows of singular vector fields introduced
in \cite{flows}.

Our plan of exposition is as follows. In section \ref{gga} we present
some fundamental results on generalized flows of distributional and
Colombeau vector fields and compare the two approaches.
Finally, section \ref{gigf} addresses the main topic of this article,
namely the analysis of group invariance under smooth Lie group actions in
the Colombeau setting.

To conclude this introduction we fix some notation and terminology
to be employed throughout the paper.  By $X$ we denote a smooth,
connected, paracompact Hausdorff manifold of dimension $n$.
For notations from Colombeau's theory of generalized functions
we follow \cite{book}. Similarly, concerning
terminology from symmetry group analysis our standard references are
\cite{Olv} for the classical theory and again \cite{book} for
the Colombeau setting.

\section{Generalized flows} \label{gga}
Lie group actions on differentiable manifolds are determined by their
one-parameter subgroups which in turn are given as the flows of
the corresponding infinitesimal generators. Therefore, in order to
analyze generalized group actions we need a theory of ordinary
differential equations with distributional or Colombeau generalized
right hand side. These problems have been addressed, e.g.,
in \cite{marsden} for the distributional and in \cite{flows}
for the Colombeau case.

We begin with a purely distributional approach,
 as provided by Marsden \cite{marsden}, and consider the initial value problem:
 \begin{equation} \label{eq:distode}
  \begin{array}{rcl}
   \dot x(t)&=&\zeta(x(t)) \\[3pt]
   x(t_0)&=& x_0,
  \end{array}
 \end{equation}
 where $\zeta \in \cd'(X,TX)$ is a distributional vector field.
We immediately note that, in the linear theory of distributions,
it is difficult to make sense of the above problem: it describes
a prospective solution $x$ which ought to be distributional, yet
take values in a manifold. Moreover, the right hand side of the
equation involves the composition of two distributional quantities.
To circumvent these conceptual problems,
Marsden's approach is to  first approximate $\zeta$ by a sequence
of smooth vector fields $\xi_\eps$. The classical theory of ODEs then
yields a flow $\Phi^\eps$ of each $\xi_\eps$. Then, under certain
assumptions on $\zeta$ and the regularizing sequence $\xi_\eps$, the
limiting measurable function $\Psi = \lim_{\eps \to 0}\Phi^\eps$
exists and is called the flow of $\zeta$.
To be more precise, let $\zeta\in\D'(X,TX)$ be a distributional
 vector field on the manifold $X$ and let $(\xi_\eps)_\eps$ be
 a net of smooth vector fields with complete flows
 $\Phi^\eps(t,.)$ and $\xi_\eps\to\zeta\in\D'(X,TX)$. $\zeta$
 is called a {\em vector field with measurable flow} $\Psi_t$ if
 \begin{itemize}
 \item [(i)] $\Phi^\eps(t,.)\to\Psi(t,.)$ almost everywhere on $X$
             for all $t$ (in particular, $\Psi_t$ is measurable), and
 \item [(ii)] For each $t\in \R$ and each $C\comp X$ there exists
   $\eps_0\in I$ and $K\comp X$ with $C\subseteq K$ such that
   $\Phi^\eps(t,C)\subseteq K$ for all $\eps$.
 \end{itemize}
It should be noted that in our terminology, (ii) says that
$\Phi^\eps(t,\,.\,)$ is c-bounded. It turns out, however,
that (in contradiction to Thm.\ 6.2 in \cite{marsden}) the assumption on
$\zeta$ to be a vector field with measurable flow $\Psi_t$
does not guarantee the flow properties of $\Psi_t$ almost
everywhere (see \cite{flows}, Prop.\ 4.1 for a counterexample).

Despite this seeming impasse, the basic approach of Marsden (i.e.,
regularizing the distributional vector field and considering the
net of flows corresponding to these regularizations), can successfully
be transferred to the Colombeau setting.

Recall that ${\gs}^1_0(X)$ denotes the space of Colombeau generalized
vector fields on $X$.
In order to derive existence and uniqueness theorems for flows of
generalized vector fields we will need the following
notions of boundedness in terms of Riemannian metrics on $X$.

 \bd Let $\xi\in\gs^1_0(X)$.
 \begin{itemize}
 \item [(i)] We say that $\xi$ is {\em locally bounded}
  resp.\ locally of {\em $L^\infty$-log-type}
  if for all $K\comp X$ and one (hence every) Riemannian metric $h$ on $X$
  we have for one (hence every) representative $\xi_\eps$
  \[
   \sup_{p\in K}\|\,\xi_\eps|_p\,\|_h\,\leq\,C\, \quad \mbox{resp.}
   \quad \sup_{p\in K}\|\,\xi_\eps|_p\,\|_h\,\le C |\log \eps|\,,
  \]
  where $\|\quad\|_h$ denotes the norm induced on $T_pX$ by $h$.
 \item [(ii)] $\xi$ is called
  {\em globally bounded} with respect to $h$ if for some (hence every) representative
  $(\xi_\eps)_\eps$ of $\xi$ there exists $C>0$ with
  \[
   \sup_{p\in X}\|\,\xi_\eps|_p\,\|_h\,\leq\, C \,.
  \]
 \end{itemize}
 \ed

Note that, contrary to local boundedness resp.\ local $L^\infty$-log type, global boundedness
obviously depends on the choice of Riemannian metric. We now have the following basic existence and uniqueness
result for ordinary differential equations with generalized right hand side on a differentiable
manifold:

 \bthm\label{th:mfode}
  Let $(X,h)$ be a complete Riemannian manifold, $\tilde x_0\in \tilde X_c$
  and $\xi\in\gs^1_0(X)$ such that
  \begin{itemize}
  \item[(i)] $\xi$ is globally bounded with respect to $h$.
  \item[(ii)] For each differential operator $P\in{\cal P}(X,TX)$ of first order $P\xi$
  is locally of $L^\infty$-log-type.
  \end{itemize}
  Then the initial value problem
 \bea\label{eq:mfode}
   \dot x(t)&=&\xi(x(t)) \\
   x(t_0)&=&\tilde x_0\nn
 \eea
  has a unique solution $x$ in $\gs[\R,X]$.
 \ethm
\pr See \cite{flows}, Th.\ 3.5. \ep
Moreover, we have:
 \bthm\label{th:mfflow}
 Let $(X,h)$ be a complete Riemannian manifold and suppose that $\xi\in\gs^1_0(X)$ satisfies
 conditions (i) and (ii) of Theorem \ref{th:mfode}.
 Then there exists a unique generalized function $\Phi\in\gs[\R\times X,X]$,
 the {\em generalized flow of $\xi$}, such that
 \bea\label{eq:mfflow1}
  \frac{d}{dt}\Phi(t,x)&=&\xi(\Phi(t,x))\quad\mbox{ in }\G^h[\R\times X,TX] \\
  \Phi(0,.)&=&\mathrm{id}_{X} \quad\mbox{ in }\G[X,X]   \label{eq:mfflow2}\\
  \Phi(t+s,.)&=&\Phi(t,\Phi(s,.)) \quad\mbox{ in }\G[\R^2\times X,X]\,.\label{eq:mfflow3}                        \eea
 \ethm
\pr See \cite{flows}, Th.\ 3.6. \ep
 \bd \label{def:flow vf}
 Let $\xi \in \cG^1_0(X)$ be a generalized vector field for which
 there exists a unique global generalized flow $\Phi \in \cG[\R\times X,X]$
 satisfying (\ref{eq:mfflow1}), (\ref{eq:mfflow2}) and (\ref{eq:mfflow3}).
 Then both $\xi$ and its generalized flow $\Phi$ are called $\cG$-complete.
$\Phi$ is called a generalized (one-parameter) group action on $X$ and
$\xi$ is called the infinitesimal generator of $\Phi$.
 \ethi

Once the existence of the generalized flow of an element $\xi$ of ${\gs}^1_0(X)$
is secured, the question arises whether there exist distributional limits
of the corresponding flows $\Phi_\eps$. For a detailed analysis of this
question (which, in a certain sense provides a resolution to the
problems encountered in the distributional modelling of generalized flows above)
we refer to \cite{flows}, sec.\ 6.

\section{Group invariants in the Colombeau setting} \label{gigf}
 If $\Phi$ is a generalized group action on $X$
 we call $u\in \gs(X)$ invariant under $\Phi$ if
 $ u(\Phi(\tilde{\eta},\tilde{x})) = u(\tilde{x})\  \forall
 \tilde{\eta} \in \widetilde{\R}_c,\ \tilde{x} \in \widetilde{X}_c$.
By the point value characterization of Colombeau generalized functions
this condition is equivalent to $u\circ\Phi = u\circ \pi_2$ as elements
of $\gs(\R\times X)$ (with $\pi_2:\R\times X \to X$ the projection).
The basic infinitesimal criterion for invariance is given in the following
proposition, proved here for the sake of completeness (cf.\ e.g., \cite{book},
Th.\ 4.5.1).
 \bp \label{prop:1}
 Let $u\in\cG(X)$ and let $\Phi$ be a generalized group
 action on $X$ with infinitesimal generator $\xi$. Then
 the following statements are equivalent:
 \begin{itemize}
 \item[(i)] $u$ is $\Phi$-invariant.
 \item[(ii)] $\xi(u)=0$ in $\cG(X)$.
 \end{itemize}
 \ethi

 \pr
 (i)$\Rightarrow$(ii) Since $u$ is $\Phi$-invariant we have
 $
 0=\frac{d}{d\eta}{\vert}_{_{0}}(u(\Phi(\eta,x)))=
 \xi(u)|_{x}$ in $\cG(X)$.

 (ii)$\Rightarrow$(i) Conversely, let $\xi(u)=0$ in $\cG(X)$. Then
 $$
 \frac{d}{d\eta}u(\Phi(\eta,\tilde{x})) = \xi(u)|_{\Phi(\eta,
 \tilde{x})} =0 \qquad \mbox{in } \gs(\R) \ \forall\tilde{x} \in
 \widetilde{X}_c.
 $$
 Therefore, for each $\tilde{x}$ the map $\eta \mapsto u(\Phi(\eta,
 \tilde{x}))$
 is constant in $\gs(\R)$, so $u\circ\Phi = u\circ \pi_2$, by
 \cite{book},
 Th.\ 3.2.8.
 \ep

To analyze the concept of group invariance for Colombeau generalized functions
let us first consider the case of a classical (i.e., smooth) generator.
If $\xi\in \mathfrak{X}(X)$ is $\gs$-complete then the generalized flow of $\xi$
coincides with the classical flow. Important examples of $\gs$-complete smooth
vector fields include:
\bexs \label{giexs}
\begin{itemize}
\item[(i)] Let $(X,g)$ be a complete Riemannian manifold and let $\xi\in \mathfrak{X}(X)$
be globally bounded w.r.t.\ $g$. Then $\xi$ is $\gs$-complete. In fact, since
$\xi_\eps \equiv \xi$ is a representative of $\xi$ as an element of ${\gs}^1_0(X)$ it is
clear that for each first order differential operator $P\in {\cal P}(X,TX)$, $P\xi$
is locally of $L^\infty$-log-type. The claim therefore follows from Th.\ \ref{th:mfflow}.
\item[(ii)] As a particular case of (i), choose $X=\R^n$ with the standard Euclidean
metric. It follows that if $\xi$ is a smooth vector field on $\R^n$ with globally
bounded coefficients then $\xi$ is $\gs$-complete and its flow is just the classical
smooth flow.
\item[(iii)] Suppose that $X=\R^n$ and the coefficients of $\xi$ are linear functions
of $x$. Then $\xi$ is $\gs$-complete (see, e.g., \cite{Ligeza2}, Th.\ 3.1). For example,
if $\xi = x_i \partial_{x_j} - x_j\partial_{x_i}$ is the smooth generator of a rotation
in $\R^n$ then it follows that $\xi$ is $\gs$-complete.
\end{itemize}
\ethi
For the base vector fields $\partial_{x_i}$ on $\R^n$, the following result by
M.\ Oberguggenberger gives a characterization of translational invariance:
\bt \label{tranlth}
Let $u\in \gs(\R^n)$. The following are equivalent:
\begin{itemize}
\item[(i)] $u(\tilde{x}_1+\eta,\tilde{x}_2,\dots,\tilde{x}_n) =
 u(\tilde{x})$ for all $\tilde{x}
\in \widetilde{\R}^n_c, \eta \in \widetilde{\R}_c$.
\item[(ii)] $\partial_{x_1} u = 0$ in $\gs(\R^n)$.
\item[(iii)] $u$ has a representative $(u_\eps)_\eps$ such that
$\partial_{x_1}u_\eps\equiv 0$
for all $\eps$.
\end{itemize}
\ethi
\pr See \cite{MoRot}, Th.\ 2.2.\ \ep
\brem \label{tranlthrem}
\begin{itemize}
\item[(i)] It follows easily from the proof of \cite{MoRot}, Th.\ 2.2.\ that
an analogous statement is valid for iterated derivatives
$\partial_{x_{i_1}}\dots \partial_{x_{i_k}}u$. 
\item[(ii)]
Until recently it was an open question whether (i)--(iii) in the above theorem is equivalent to
\begin{itemize}
\item[(i')] $u(\tilde{x}_1+\eta,\tilde{x}_2,\dots,\tilde{x}_n) =
 u(\tilde{x})$ for all $\tilde{x}
\in \widetilde{\R}^n_c, \eta \in \R$.
\end{itemize}
i.e., whether standard translations suffice to characterize
translational invariance of Colombeau generalized functions.
In \cite{PSV}, however, Pilipovi\'c, Scarpalezos and
Valmorin were able to give an affirmative answer to this question.
\end{itemize}
\ethi
Note that Th.\ \ref{tranlth} implies in particular that a generalized function $u$ is invariant
under translations if and only if it possesses a distinguished representative $(u_\eps)_\eps$
such that each $u_\eps$ is a translation invariant smooth function. The following proposition locally extends the
validity of this result to all smooth group actions which are regular in the
sense that their infinitesimal generators are non-vanishing.
\bp \label{gentransprop}
Let $\xi\in \mathfrak{X}(X)$ be $\gs$-complete with flow $\Phi$ and suppose that $\xi(x_0)\not=0$ for
some $x_0\in X$. Then there exists a neighborhood $U$ of $x_0$ in $X$ such that $u\in \gs(U)$ is invariant
under $\Phi$ if and only if $u$ possesses a representative $(u_\eps)_\eps$ with each $u_\eps$
invariant under $\Phi$.
\ethi
\pr By \cite{AM}, Th.\ 2.1.9 we may ``straighten out'' $\xi$ around $x_0$. Thus there exists a chart
$(U,\vphi)$ around $x_0$ with $\vphi_*\xi = \partial_{x_1}$. But then $\vphi_* u|_U$ and hence
$u|_U$ itself satisfy the claim by Th.\ \ref{tranlth}.
\ep
We next wish to generalize Prop.\ \ref{gentransprop} from one-parameter groups
to more general Lie group actions on manifolds. Thus let us assume that $G$ is a Lie group and
$\Psi: G\times X \supseteq \mathcal{V} \to X$ is a regular transformation group
(i.e., all orbits have the same dimension as submanifolds and each point in $X$ has
a base of neighborhoods whose elements intersect each orbit in a connected subset thereof,
cf.\ \cite{Olv2}, p.\ 41). Analogous to the smooth setting, we call an element $u$ of
$\gs(X)$ invariant under $\Psi$ if
$$
u\circ \Psi = u\circ \pi_2 \qquad \mbox{in } \gs(\mathcal{V}) \,,
$$
where $\pi_2: G\times X \to X$ is the projection onto the second factor. The desired
generalization of Prop.\ \ref{gentransprop} then takes the following form:
\bt \label{gengroup}
Let $\Psi: G\times X \supseteq \mathcal{V} \to X$ be a regular Lie group action
on $X$ with $s$-dimensional orbits. Then each point $x_0\in X$ possesses a neighborhood $U$
such that the following statements are equivalent for each $u\in \gs(U)$:
\begin{itemize}
\item[(i)] $u$ is invariant under $\Psi$.
\item[(ii)] $u$ possesses a representative $(u_\eps)_\eps$ with each $u_\eps$
invariant under $\Psi$.
\end{itemize}
\ethi
\pr By \cite{Olv2}, Th.\ 2.23 we may choose a rectifying local chart
$\vphi: x\mapsto (y_1,\dots,y_s,z_1,\dots,z_{n-s})$ in a neighborhood $U$ of
$x_0$ such that each group orbit intersects
the coordinate chart in at most one slice $\{z_1=c_1,\dots,z_{n-s}$ $=$ $c_{n-s}\}$
(with $c_i$ constant for $1\le i \le n-s$).  Denote by $O_x$ the orbit of
$\Psi$ through $x$.
Then for each $x\in U$ the
local vector fields $\left. \frac{\partial}{\partial y_i}\right|_x$
($1\le i\le s$) form a basis of the tangent space of $O_x$.
It follows that in this neighborhood each infinitesimal generator of $\Psi$
is a unique $\cinfty$-linear combinations of the
$\frac{\partial}{\partial y_i}$ and vice versa. In $U$, invariance under $\Psi$
therefore amounts to $\partial_{y_i}u$ being zero in $\gs(U)$. Therefore, by an
application of Remark \ref{tranlthrem} (i) to $\vphi_*u$ we reach the desired
conclusion.
\ep
As an important concrete example of invariance of Colombeau generalized functions
under smooth transformation groups let us consider in some detail the
case of rotational invariance, following \cite{ObCont,MoRot,KoK}.
Let $\mathsf{SO}(n,\widetilde{\R})$ denote the special orthogonal group over the
ring $\widetilde{\R}$ of generalized numbers and $\mathsf{SO}(n,\R)$ the
usual special orthogonal group. Letting $\mathsf{SO}(n,\widetilde{\R})$ act
naturally on $\R^n$ it is well known that a basis of the Lie algebra of infinitesimal
generators of this action is given by the set of all $\xi_{ij} = x_i\partial_{x_j}
- x_j\partial_{x_i}$ for $i<j$. For an element $u$ of $\R^n$ to be invariant under
the flow $\Phi_{ij}$ of $\xi_{ij}$ means that
$$
u(\Phi(\tilde \eta,\tilde x)) = u(\tilde x) \qquad \forall \tilde \eta \in \tilde \R_c \ \forall \tilde x \in \tilde \R_c^n\,.
$$
Here, the action of $\Phi_{ij}(\tilde \eta,\,.\,)$ is precisely a rotation in the $(x_i,x_j)$-plane by the
generalized angle $\tilde \eta$, hence is given by the action of the corresponding element
of $\mathsf{SO}(n,\widetilde{\R})$ on $\tilde x$. Conversely, as was shown in \cite{MoRot}, Sec.\ 2, Lemma 3, each
generalized rotation $A \in \mathsf{SO}(n,\widetilde{\R})$ in the $(x_i,x_j)$-plane
is precisely of this form. (This structural relationship in fact reaches even further:
by the same result of Oberguggenberger, the $\xi_{ij}$ also form a basis of the ``Lie algebra''
of $\mathsf{SO}(n,\widetilde{\R})$ in the following sense: each $A\in \mathsf{SO}(n,\widetilde{\R})$
is of the form $\exp(v)$ for some generalized vector field $v=\sum_{i<j}\alpha_{ij}\xi_{ij}$ with
$\alpha_{ij}\in \tilde\R$ for all $i<j$.)
Consequently, the action of a generalized $(i,j)$-rotation $A\in \mathsf{SO}(n,\widetilde{\R})$
(which is an example of a generalized group action in the sense of section \ref{gga})
can be viewed as the ``nonstandardization'' of the corresponding classical rotation which is obtained
by replacing the real angle $\eta$ by the generalized angle $\tilde \eta$. This, of course, is
a direct result of the $\gs$-completeness of the smooth generators $\xi_{ij}$ (cf.\ Ex.\ \ref{giexs} (iii)).
Combining these observations with Th.\ \ref{prop:1} we obtain
(see \cite{MoRot} for an alternative direct proof):
\bp \label{rotprop}
Let $u\in \gs(\R^n)$. The following are equivalent:
\begin{itemize}
\item[(i)] $u\circ A = u$ for all $A\in \mathsf{SO}(n,\widetilde{\R})$, i.e.,
$u$ is rotationally invariant.
\item[(ii)] $\xi(u) = 0$ for each infinitesimal generator $\xi$ of $\,\mathsf{SO}(n,\R)$.
\end{itemize}
\ethi
Moreover, we have:
\bp For $u\in \gs(\R^n\setminus \{0\})$, (i) and (ii) are further equivalent with
\begin{itemize}
\item[(iii)] $u$ possesses a representative consisting entirely of rotation invariant functions.
\end{itemize}
\ethi
\pr $\mathsf{SO}(n,\widetilde{\R})$ acts freely on $\R^n\setminus \{0\}$, so we may
employ Prop.\ \ref{gentransprop} to establish the claim.
\ep
The restriction to $\R^n\setminus \{0\}$ in the above result is grounded in the method of proof (application of
Prop.\ \ref{gentransprop}) rather than in the subject matter itself. In fact, the equivalence is
true on all of $\R^n$ (see \cite{MoRot}).

The above chain of equivalences raises the question whether in the present context the
analogue of Remark \ref{tranlthrem} (ii) can be established as well, i.e., whether (i)--(iii)
are equivalent to
\begin{itemize}
\item[(i')]  $u\circ A = u$ for all $A\in \mathsf{SO}(n,\R)$.
\end{itemize}
In fact, we have: \bp \label{rn0prop} For $u\in
\gs(\R^n\setminus\{0\})$, (i) and (i') are equivalent. \ethi \pr
Since each $\xi_{ij}$ is nonzero on $\R^n\setminus \{0\}$ it may
be straightened out, in fact even globally on all of
$\R^n\setminus \{0\}$ (by using appropriate polar coordinates as
charts), cf.\ the proof of Prop.\ \ref{gentransprop}. This
procedure reduces the proof to an application of Remark
\ref{tranlthrem}. \ep A direct extension of  the proof of Prop.\
\ref{rn0prop} to the case $X=\R^n$ is not possible: contrary to
the smooth situation a Colombeau generalized function which is
rotationally invariant on $\R^n\setminus \{0\}$ need not be
rotationally invariant on all of $\R^n$. As an example, choose any
test function $\vphi$ whose support is not rotationally invariant
and set $u=[(\vphi(\,.\,/\eps))_\eps]$. Then $u$ is supported in
$\{0\}$ yet it is clearly not rotationally invariant on $\R^n$.
Despite this technical complication, however, it turns out that
the {\em result} can be extended to all of $\R^n$, thereby
providing an affirmative answer to a question raised by M.\
Oberguggenberger in \cite{MoRot}: \bt For $u\in \gs(\R^n)$, (i)
and (i') are equivalent. \ethi \pr This result was established in
\cite{KoK}, Th.\ 7.6. We include a proof here for the reader's
convenience. It clearly suffices to show that (i') implies (i).
Let us first consider the case $n=2$. Let $\widetilde{A}\in
\mathsf{SO}(2,\widetilde{\R})$. Then by the discussion preceding
Prop.\ \ref{rotprop} there exists some $\tilde{\eta} \in
\widetilde{\R}_c$ such that
$$
\widetilde{A} = \left[ \left(
\begin{array}{cr}
\cos(\eta_\eps) & -\sin(\eta_\eps) \\
\sin(\eta_\eps) & \cos(\eta_\eps)
\end{array}
\right)_{\!\!\eps}\right]\,.
$$
Given $\tilde{x}$, $\tilde{y} \in \widetilde{\R}_c$ we have to
show that $u(\widetilde{A} \cdot (\tilde{x},\tilde{y})^t) =
u(\tilde{x},\tilde{y})$ in $\widetilde{\R}$. We may write
$(\tilde{x},\tilde{y}) = [(r_\eps
\cos(\theta_\eps),r_\eps\sin(\theta_\eps))]$ for suitable $r_\eps
\ge 0$, $\theta_\eps$. Now set $v_\eps := \theta \mapsto
u_\eps(r_\eps \cos(\theta),r_\eps\sin(\theta))$. Then
$v=[(v_\eps)_\eps]\in \gs(\R)$ and by assumption $v(\tilde{\theta}
+\eta) = v(\tilde{\theta})$ in $\widetilde{\R}$ for all
$\tilde{\theta} \in \widetilde{\R}_c$ and all $\eta\in \R$. But
then by Th.\ \ref{tranlth} it follows that $v$ is a generalized
constant, thereby finishing the proof for $n=2$.

In the general case $n\ge 2$ we verify (ii) of Prop.\ \ref{rotprop}. Let
$1\le i < j \le n$ and let $ \xi_{ij}  = x_i\partial_{x_j} - x_j
\partial_{x_i}$ as above be
an infinitesimal generator of $\mathsf{SO}(n,\R)$. Fix compactly
supported generalized numbers $\tilde{x}_1$,\dots,$\tilde{x}_{i-1}$,
$\tilde{x}_{i+1}$,\dots, $\tilde{x}_{j-1}$,
 $\tilde{x}_{j+1}$,\dots,$\tilde{x}_{n}$ and consider the maps
$$
w_\eps: (x_i,x_j) \mapsto u_\eps(\tilde{x}_1,\dots,
\tilde{x}_{i-1},x_i,\dots,\tilde{x}_{j-1},x_j,\dots,\tilde{x}_n)\,.
$$
Then $w=[(w_\eps)_\eps] \in \gs(\R^2)$ and from our assumption it
follows that $w\circ A = w$ in $\gs(\R^2)$ for all
$A\in \mathsf{SO}(\R^2)$. By what we have already proved in the 2D-case
and Prop.\ \ref{rotprop} it follows that
$ \xi_{ij}  w = 0$ in $\gs(\R^2)$ for each $i<j$, which finishes the proof.
\ep

In the smooth setting, the local structure of invariants of a group action is
determined by so-called complete sets of functionally independent invariants
(cf.\ \cite{Olv}, Sec.\ 2.1). A family of smooth functions $f_1,\dots,f_k$ on $X$
is called functionally dependent if for each $x\in X$ there exists a neighborhood
$U$ of $x$ and a smooth function $F$ which is not identically zero on any open
subset of $\R^k$ such that $x\mapsto F(f_1(x),\dots,f_k(x))$ vanishes identically
on $U$. $f_1,\dots,f_k$ are called functionally independent if they are not
functionally dependent on any open subset of $X$. It is shown in \cite{Olv},
Th.\ 2.17 that if $\Psi$ acts (semi-)regularly with $s$-dimensional orbits
then in a neighborhood of any point $x_0$ of $X$ there exists a set
$f_1,\dots,f_{n-s}$ of functionally independent invariants such that
any other local invariant $f$ of $\Psi$ is of the form
$$
f(x) = F(f_1(x),\dots,f_{n-s}(x))
$$
for some smooth function $F$. Such a family $f_1,\dots,f_{n-s}$ is called
a complete set of functionally independent invariants of $\Psi$.
Using Th.\ \ref{gengroup} we now show that an analogous characterization
of generalized invariants of regular smooth group actions holds true.

\bt \label{fiinv}
Let $\Psi: G\times X \supseteq \mathcal{V} \to X$ be a regular Lie group action
on $X$ with $s$-dimensional orbits. Then for each $x_0\in X$ there exists a
neighborhood $U$ of $x_0$ and a complete set of functionally independent
invariants $f_1,\dots,f_{n-s}$ on $U$ such that every $\Psi$-invariant
$u\in \gs(U)$ is of the form $u(x)=v(f_1(x),\dots,f_{n-s}(x))$ for some
$v \in \gs(V)$, where $V$ is an appropriate open subset of $\R^{n-s}$
\ethi
\pr Employing the notations of Th.\ \ref{gengroup}, the coordinates $z_1=f_1(x)$,
\dots $z_{n-s}=f_{n-s}(x)$ are in fact a complete set of functionally independent invariants
for the group action $\Psi$ on $X$
(cf.\ the proof of Th.\ 2.17 in \cite{Olv}). Now if $u\in \gs(U)$ is an invariant
of $\Psi$ then by Th.\ \ref{gengroup} we may choose a representative $(u_\eps)_\eps$
of $u$ such that each $u_\eps$ is an invariant of $\Psi$ on $U$.
Then by \cite{Olv},
Th.\ 2.17, for each $\eps$ there exists a smooth function $v_\eps$ such that
$u_\eps(x) = v_\eps(f_1(x),\dots,f_{n-s}(x))$. Denoting $\vphi_*u$ by $\tilde u$
we have $\tilde u_\eps(y,z) = v_\eps(z)$, so $(v_\eps)_\eps \in \esm(V)$,
where $V$ is the projection of $\vphi(U)$ onto the last $n-s$ components.
It follows that the class $v=[(v_\eps)_\eps]$ of $(v_\eps)_\eps$ in $\gs(V)$ is the
desired generalized function.
\ep


 \end{document}